\documentclass[12pt,psfig]{article}
\usepackage{bibentry}
\usepackage{amsfonts}
\usepackage{amssymb}
\usepackage{amsthm}
\usepackage[namelimits]{amsmath}
\usepackage{graphicx,epsfig}
\usepackage{amscd}
\usepackage{epsfig,color,psfrag}
\date{}

\newcommand{\tss}{\textstyle}

\newcommand{\ep}{\epsilon}

\newcommand{\ga}{\gamma}

\newcommand{\al}{\alpha}
\newcommand{\be}{\beta}

\newcommand{\si}{\sigma}
\newcommand{\ld}{\lambda}
\newcommand{\Ld}{\Lambda}
\newcommand{\ba}{\begin{array}}
\newcommand{\ea}{\end{array}}
\newcommand{\beq}{\begin{eqnarray*}}
\newcommand{\eeq}{\end{eqnarray*}}
\newcommand{\bdm}{\begin{displaymath}}
\newcommand{\edm}{\end{displaymath}}

\theoremstyle{definition}
\newtheorem{theorem}{Theorem}[section]

\newtheorem{lemma}{Lemma}[section]

\newtheorem{remark}{Remark}[section]

\begin{document}

\pagestyle{plain}

\begin{center}
{\large \bf  Hyperbolicity of the invariant sets for the real polynomial maps}
\\ [0.3in]

XU ZHANG \footnote{ Email address:\ \ xuzhang08@gmail.com (X. Zhang).}

 \vspace{0.15in}
{\it  Department of Mathematics, Shandong University \\[-0.5ex]
Jinan, Shandong 250100, P.~R. China}

\end{center}

\vspace{0.18in}

\baselineskip=20pt

{\bf\large{ Abstract.}}\ It is well known that for $a>4$, the dynamical behaviors of the logistic map $f_a(x)=ax(1-x)$ on the maximal invariant compact set are ``simple", which could be clearly explained by the theories of hyperbolic dynamics and symbolic dynamics. Is it possible that similar phenomena could be observed in general real polynomial maps? In this paper, we study this problem by investigating the real polynomial map $f_a(x)=ag(x)$, where $a$ is a parameter, and $g$ is a real-coefficient polynomial, which has at least two different real zeros or only one real zero.
\bigskip

\bigskip


\section{Introduction}

Many books, introduction to dynamical systems, would take an account of the logistic map $f_a(x)=ax(1-x)$. For any $a>4$, there exists a maximal compact invariant set $\Ld_a$, on which the map is hyperbolic and topologically conjugate to the one-sided fullshift on two symbols. When $a>2+\sqrt{5}$, it is easy to obtain this conclusion, since the absolute value of the derivative of every point in $\Ld_a$ is greater than $1$, however, for $4<a\leq2+\sqrt{5}$, it is a little hard \cite{Robinson}.

The study of the logistic map with $a>4$ dates back to the work of Fatou and Julia on complex dynamics, who obtained that the Julia set of a polynomial map would be totally disconnected if all the critical points go to infinity under the iteration of the map \cite{CarlesonGamelin}. By using method in terms of real variables, Henry proved that almost every point in the unit interval would escape from it under the iteration of the logistic map \cite{Henry}. Guckenheimer \cite{Guckenheimer}, Misiurewicz \cite{Misiurewicz} and van Strien \cite{Strien} had a proof with the fact that the logistic map has negative ``Schwarzian derivative". In \cite{Robinson}, Robinson gave an elegant proof with the help of the Schwarz Lemma. Recently, Aulbach and Kieninger \cite{Aulbach},
Glendinning \cite{Glendinning} and Kraft \cite{Kraft} have given different elementary proofs on the hyperbolicity of the logistic map with $a>4$.

However, the general polynomial maps have received less attention. Motivated by these important works on the logistic map, we study a general real polynomial map $f_a(x)=a(x-b_1)^k(b_2-x)^lh(x)$, where $b_1$ and $b_2$ are constants with $b_1\neq b_2$; $k$ and $l$ are positive integers; $a$ is a parameter; and $h$ is a polynomial satisfying that $|h(x)|>0$ for $x\in[\min\{b_1,b_2\},\max\{b_1,b_2\}]$. We first investigate the map under the assumptions that $b_1b_2\geq0$, $|b_1|<|b_2|$, and $a(\mbox{sign}(b_2))^{l+1}h(b_1)>0$. Except the situation that $b_1=0$ and $l>k=1$, we show that for sufficiently large $|a|$, the polynomial map $f_a$ has an invariant Cantor set on which it is hyperbolic and topologically conjugate to the one-sided fullshift on two symbols (see Theorems \ref{2positive} and \ref{zero}). For the polynomial map $f_a$, when $|a|$ is large enough, one could find two disjoint compact intervals $I_1$ and $I_2$ dependent on $a$, such that $[\min\{b_1,b_2\},\max\{b_1,b_2\}]\supset f_a(I_1)=f_a(I_2)\supset I_1\cup I_2$. In the three cases that $\rm(1)$ $b_1b_2>0$, $\rm(2)$ $b_1b_2=0$ and $k>l\geq1$, $\rm(3)$ $b_1b_2=0$ and $k=l=1$, we obtain that there exists a constant $\ld>1$, such that $\min_{x\in I_1\cup I_2}|f'_a(x)|\geq\ld$ for sufficiently large $|a|$ by applying the elementary methods. However, in the case that $b_1b_2=0$ and $l>k\geq2$, we find that $|f'_a(x)|>0$ for all $x\in I_1\cup I_2$ and $\min_{x\in I_1\cup I_2}|f'_a(x)|\to0$ as $|a|\to\infty$, which is different from what we have observed in the logistic map, since the logistic map satisfies the assumption that $b_1=0$, $b_2=1$, $k=l=1$, and $h(x)\equiv1$. To overcome the problem that we could not obtain good estimation on the derivative of $f_a(x)$ on $I_1\cup I_2$ in the case that $b_1=0$ and $l\geq k\geq2$, we utilize the methods in complex dynamics to prove the hyperbolicity of the invariant sets. Moreover, in the case that $b_1=0$ and $l>k=1$, we show that when $|a|$ is large enough, there exists an invariant set on which $f_a$ is topologically semiconjugate to the one-sided fullshift on two symbols but $f_a$ is not hyperbolic on it (see Theorem \ref{nonhyp}). Then we study the map under the assumptions that $b_1b_2\geq0$, $|b_1|<|b_2|$, and $a(\mbox{sign}(b_2))^{l+1}h(b_1)<0$. We obtain that there exists a hyperbolic invariant set for $f_a$ on which $f_a$ is topologically conjugate to the one-sided fullshift on two symbols under certain conditions (see Theorems \ref{2positiveneg} and \ref{zeropositive}). Finally, in the case that $b_1b_2<0$, we show that there exists a hyperbolic invariant set for $f_a$ when $|a|$ is large enough (see Theorem \ref{diffroot}). And, for the polynomial map $f_a(x)=ag(x)$, where $g$ has only one real zero, we show that there exists a hyperbolic invariant set for $f_a$ under certain conditions (see Theorem \ref{single}).

The rest of the paper is organized as follows. Section 2 contains
some basic concepts and useful results. In Section 3, the polynomial map $f_a(x)=ag(x)$ is investigated, where $g$ has two distinct non-negative or nonpositive real zeros. In Section 4, the polynomial map $f_a(x)=ag(x)$ is studied, where $g$ has one positive and one negative real zeros, or only one real zero. \bigskip

\section{Preliminaries}

In this section, basic definitions and useful results are introduced. \medskip

We first introduce several notions about matrix and the definition of the one-sided symbolic dynamical system \cite{Robinson}. A matrix $A=(a_{ij})_{m\times
m}$ $(m\geq2)$ is said to be a transition matrix if $a_{ij}=0$
or $1$ for all $i, j;\ \sum^m_{j=1}a_{ij}\geq1$ for all $i$; and
$\sum^m_{i=1}a_{ij}\geq1$ for all $j$, $1\leq i,j\leq m.$
$A$ is called positive, if all its entries $a_{ij}>0$.
$A$ is said to be eventually positive if there exists a positive
integer $k$ such that $A^n>0$ for all the integers $n\geq k$.

Let $S_0:=\{1,2,...,m\},\ m\geq2$, and $\textstyle\sum_m:=\{\al=(a_0,a_1,a_2,...):\ a_i\in S_0,\ i\geq0\}$ be the one-sided sequence space. We define a metric
on $\textstyle\sum_m$ by
$$ d(\al,\beta)=\sum^{\infty}_{i=0}\frac{\delta(a_i,b_i)}{2^i},$$
where $\al=(a_0,a_1,a_2,...)$, $\beta=(b_0,b_1,b_2,...)\in\sum_m$,
$\delta(a_i,b_i)=1$ if $a_i\neq b_i$, and $\delta(a_i,b_i)=0$ if $a_i=b_i$,
$i\geq0$. Then, $(\sum_m,d)$ is a complete metric space. Define the shift map $\sigma:\sum_m\to\sum_m$ by
$\sigma(\al)=(a_1,a_2,...)$, where $\al=(a_0,a_1,...)$.
Then, $(\sum_m,\sigma)$ is called the one-sided symbolic dynamical
system on $m$ symbols.

Given a transition matrix $A=(a_{ij})_{m\times m}$, denote
$$\textstyle\sum_A:=\{\be=(b_0,b_1,...)\in\sum_m:a_{b_ib_{i+1}}=1,\ i\geq0\}.$$
The map
$$\si_A:=\si|_{\sum_A}:\textstyle\sum_A\to\sum_A$$
is said to be the subshift of finite type for $A$.
Obviously, $\sum_A$ is invariant
under $\si_A$ and a compact subset of $\sum_m$. For the transition matrix $A$, a
finite sequence $w=(s_1,s_2,...,s_k)$ is said to be an allowable word of
length $k$ for $A$ if $a_{s_is_{i+1}}=1$, $1\leq i\leq k-1$, where
$s_1,s_2,...,s_k\in S_0$.

Next, the notion of hyperbolicity of one dimensional dynamics is given \cite{Kraft}. Let $f:\mathbf{R}\to\mathbf{R}$ be a $C^1$ map, and $\Ld$ be a compact invariant set for $f$, i.e., $f(\Ld)=\Ld$. Then $\Ld$ is a hyperbolic set for $f$ if there are constants $C>0$ and $\ld>1$ such that $|(f^n)'(x)|\geq C\ld^n$ for all $x\in\Ld$ and all $n\geq1$.

The next result is introduced from \cite[Lemma 4]{Kraft} and \cite[Lemma 2.1]{Aulbach}.

\begin{lemma} \label{hyp}
Let $f:\mathbf{R}\to\mathbf{R}$ be a $C^1$-map, and let $\Ld\subset\mathbf{R}$ be compact invariant set for $f$. Then $\Ld$ is hyperbolic with respect to $f$, if and only if for each $x\subset\Ld$ there exists a positive integer $k_x$ (depending on $x$) such that $|(f^{k_x})'(x)|>1$.
\end{lemma}

Finally, we introduce needed notions and useful results about complex dynamics \cite{CarlesonGamelin}.

Set $\overline{\mathbf{C}}:=\mathbf{C}\cup\{\infty\}$. Denote $|z|$ for the modulus of $z\in\mathbf{C}$. For an analytic function $p$ defined on $\overline{\mathbf{C}}$, if $p(z_0)=z_0$, then $z_0$ is called a fixed point of $p$. The number $p'(z_0)$ is said to be the multiplier of $p$ at $z_0$. If $|p'(z_0)|<1$, then $z_0$ is called an attracting fixed point; if $|p'(z_0)|=0$, then $z_0$ is a superattracting fixed point. The point $z_0$ is said to be periodic if $z_0=z_n$ for some non-negative integer $n$, where $z_n=p^n(z_0)$. The minimal $n$ is its period, the orbit $\{z_1,z_2,...,z_n\}$ is called a cycle, and the cycle is said to be attracting if $|(p^n)'(z_0)|<1$.

Let $\mathcal{F}$ be a family of meromorphic functions in a domain $D\subset\overline{\mathbf{C}}$. If every sequence $\{f_n\}$ in $\mathcal{F}$ contains a subsequence that converges uniformly in the spherical metric on compact subsets of $D$, then $\mathcal{F}$ is called a normal family. Let $R=P/Q$ be a rational map, where $P$ and $Q$ are polynomials with no common factors and degree $d=\max\{\mbox{deg}P,\ \mbox{deg}Q\}\geq2$. The Fatou set of $R$ is defined to be the set of points $z_0\in\overline{\mathbf{C}}$ such that $\{R^n\}$ is a normal family in some neighborhood of $z_0$. The complement of the Fatou set is said to be the Julia set. A set is called completely invariant for $R$ if both it and its complement are invariant. The Julia set is completely invariant \cite[Page 56]{CarlesonGamelin}. If $U$ is a completely invariant component of the Fatou set, then the boundary of $U$ is equal to the Julia set, and there are at most two completely invariant components of the Fatou set \cite[Page 70]{CarlesonGamelin}. Denote $A(z_0,R)$ as the basin of attraction of an attracting fixed point $z_0$, that is, $A(z_0,R)$ consists of $z$ such that $R^n(z)\to z_0$ as $n\to+\infty$. If $z_0$ is an attracting fixed point, then the basin of attraction $A(z_0,R)$ is a union of components of the Fatou set, and the boundary of $A(z_0,R)$ coincides with the Julia set \cite[Page 58]{CarlesonGamelin}. If $R=P$ is a polynomial of degree $d\geq2$, then $\infty$ is a superattracting fixed point of $R$ and the Julia set coincides with the boundary of $A(\infty,R)$ \cite[Page 65]{CarlesonGamelin}. The critical point of $R$ is a point on the sphere where $R$ is not locally one-to-one. The set of critical points consist of solutions of $R'(z)=0$ and of poles of $R$ of order two or higher. The rational function $R$ is hyperbolic on the Julia set if and only if every critical point belongs to the Fatou set and is attracted to an attracting cycle \cite[Page 90]{CarlesonGamelin}.

\bigskip

\section{Polynomial with two different non-negative or nonpositive real zeros}

In this section, a type of real-coefficient polynomial maps is studied, where the polynomial has at least two different non-negative real zeros. By applying the methods used in this section, similar results can be obtained when the polynomial has at least two distinct nonpositive real zeros.

Consider a real-coefficient polynomial of degree $n$,
$$g(x)=x^{n}+a_{n-1}x^{n-1}+\cdots+a_0,                    $$
which has at least two different zeros. Rewrite $g$ as
$$g(x)=(x-\al_1)^{m_1}\cdots(x-\al_{r})^{m_{r}}
(x^2+\be_{r+1}x+\ga_{r+1})^{m_{r+1}}\cdots
(x^2+\be_{s}x+\ga_{s})^{m_{s}},
$$
where
$$m_i\geq1,\ 1\leq i\leq s;\ r\geq2;\ \sum_{i=1}^r m_i+
\sum_{i=r+1}^{s} 2m_{i}=n;\ \be_{j}^2-4\ga_{j}<0,\ r+1\leq j\leq s,
$$ and $\al_1,\dots,\al_{r}$ are real zeros of $g(x)$ with
$\al_1<\cdots<\al_{r}$.

Consider the following identity:
$$\frac{\tss f_a'(x)}{\tss f_a(x)}=\frac{\tss g'(x)}{\tss
g(x)}=\sum_{i=1}^{r}\frac{\tss
m_{i}}{x-\al_{i}}+\sum_{i=r+1}^{s}\frac{\tss m_{i}(2x+\be_{i})}{\tss
x^2+\be_{i}x+\ga_{i}}.
\eqno{(3.1)}$$ It can be easily shown that
$$
\lim_{x\to\al_{i}^+}\frac{\tss f_a'(x)}{\tss
f_a(x)}=\lim_{x\to\al_{i}^+}\frac{\tss g'(x)}{\tss
g(x)}=+\infty,\;\; \lim_{x\to\al_{i}^-}\frac{\tss f_a'(x)}{\tss
f_a(x)}=\lim_{x\to\al_{i}^-}\frac{\tss g'(x)}{\tss
g(x)}=-\infty,\ \ 1\leq i\leq r.                          \eqno{(3.2)}$$
If there exists $j$, $1\leq j\leq r$, such that $m_{j}\geq2$, then $\al_{j}$ is also a real root of $g'(x)$, and one has
$$
\lim_{x\to\al_{j}^+}\frac{\tss f_a''(x)}{\tss
f_a'(x)}=\lim_{x\to\al_{j}^+}\frac{\tss g''(x)}{\tss
g'(x)}=+\infty,\;\; \lim_{x\to\al_{j}^-}\frac{\tss f_a''(x)}{\tss
f_a'(x)}=\lim_{x\to\al_{j}^-}\frac{\tss g''(x)}{\tss
g'(x)}=-\infty.                          \eqno{(3.3)}$$
In the following discussions, (3.1), (3.2) and (3.3) will play an important role.

In this section, suppose that there exists $i_0$, $1\leq i_0< r$, such that $\al_{i_0}\geq0$. We investigate the following polynomial map:
$$f_a(x)=ag(x),                  $$
where $a\in \mathbf{R}$ is a parameter.

We will first study the dynamical behaviors of $f_a$, where $a$ satisfies
$$f^{(m_{i_0})}_a(\al_{i_0})=ag^{(m_{i_0})}(\al_{i_0})>0.            $$
This assumption guarantees that $f_a(x)>0$ for all $x\in(\al_{i_0},\al_{i_0+1})$. Since  $f_a(x)$ can be written as
$$f_a(x)=\frac{\tss 1}{\tss
m_{i_0}!}f^{(m_{i_0})}_a(\al_{i_0})(x-\al_{i_0})^{m_{i_0}}+
o(|x-\al_{i_0}|^{m_{i_0}}),                                 $$
there exists no real zero of $f_a(x)$ in $(\al_{i_0},\al_{i_0+1})$ and
$f^{(m_{i_0})}_a(\al_{i_0})>0$, one has that $f_a(x)$ is positive in
$(\al_{i_0},\al_{i_0+1})$.

\begin{lemma} \label{deriva} Suppose that $\al_{i_0}>0$ and $f^{(m_{i_0})}_a(\al_{i_0})>0$. For the polynomial $f_a(x)$, there exist two points $x_0$ and $
x_0'$ with $\al_{i_0}<x_0<x_0'<\al_{i_0+1}$, such that $f'_a(x)>0$ for all $ x\in (\al_{i_0},x_0]$, and $f'_a(x)<0$ for all $ x\in [x_0',\al_{i_0+1})$. Further, if $m_{i_0}\geq2$, then there exists $x_1$, $x_1\in(\al_{i_0},x_0]$, such that $f''_a(x)>0$ for all $x\in(\al_{i_0},x_1]$;
if $m_{i_0+1}\geq2$, then there exists $x'_1$, $x'_1\in[x'_0,\al_{i_0+1})$, such that $f''_a(x)>0$ for all $x\in[x'_1,\al_{i_0+1})$.
\end{lemma}

\begin{proof}
The existence of $x_0$ and $x_0'$ can be derived from (3.1) and (3.2). This, together with (3.3), implies the existence of $x_1$ and $x'_1$.
\end{proof}
\medskip

\begin{lemma} \label{singleroot}
For the polynomial $f_a(x)$, suppose that $g'(\al_{j})\neq0$ for some $j$, $1\leq j\leq r$. Then, there exists $\delta>0$ such that $g'(x)\neq0$ for all $x\in[\al_j-\delta,\al_j+\delta]$. And, for any given constant $\ld>1$, there exists a constant $N_{\ld}>0$ such that for any $|a|>N_{\ld}$, one has $|f'_a(x)|\geq\ld$ for all $x\in[\al_j-\delta,\al_j+\delta]$.
\end{lemma}

\begin{proof}
The proof of this lemma is simple. So, it is omitted here.
\end{proof}

We come now to one of the principal theorems.

\begin{theorem} \label{2positive}
Suppose that $\al_{i_0}>0$ and $f^{(m_{i_0})}_a(\al_{i_0})>0$. For sufficiently large $|a|$, there exists an invariant Cantor set $\Ld_a=\{x:\ f^k_a(x)\in[\al_{i_0},\al_{i_0+1}]\ \mbox{for all}\ k\geq0\}$ such that $f_a$ is hyperbolic on $\Ld_a$ and $f_a:\Ld_a\to\Ld_a$ is topologically conjugate to $\si:\sum_2\to\sum_2$.

\end{theorem}

\begin{proof}
The whole proof is divided into two steps.

\indent \bf{Step}\ \rm{1}. The properties of the polynomial $f_a(x)$ are discussed when $|a|$ is large enough.

By Lemma \ref{deriva}, there exist $x'_l$ and $x'_r$ with $\al_{i_0}<x'_l<x'_r<\al_{i_0+1}$, such that $f'_a(x)>0$ for all $x\in(\al_{i_0},x'_l]$, and $f'_a(x)<0$ for all $x\in[x'_r,\al_{i_0+1})$. Choose $\ld>1$. By (3.2), there exist $x_l$ and $x_r$, $x_l\in(\al_{i_0},x'_l]$, $x_r\in[x'_r,\al_{i_0+1})$, such that
$$\bigg|\frac{ f_a'(x)}{
f_a(x)}\bigg|=\bigg|\frac{g'(x)}{g(x)}\bigg|\geq\frac{\ld}{\al_{i_0}},\ \ x\in(\al_{i_0},x_l]\cup[x_r,\al_{i_0+1}).         \eqno{(3.4)}$$

Denote
$$m_0:=\inf_{x\in[x_l,x_r]}|g(x)|.$$

Next, it is to show that if
$$|a|\geq\frac{\al_{i_0+1}}{m_0}\ \mbox{and}\ ag^{(m_{i_0})}(\al_{i_0})>0,\eqno(3.5)$$
then the assertion of this theorem holds.

Now, it is to prove that there exists one and only one fixed point of $f_a$ in $(\al_{i_0},x_l)$ for any $a$ satisfying (3.5). Obviously, $f_a(x_l)\geq\al_{i_0+1}>x_l$. This, together with $f_a(\al_{i_0})=0$, implies that there exists at least one fixed point of $f_a$ in $(\al_{i_0},x_l)$ for any $a$ satisfying (3.5). Set $z_a:=\min_{x\in(\al_{i_0},x_l)}\{x:\,f_a(x)=x\}$. It is evident that $z_a>\al_{i_0}$. Since $f'_a(x)>0$ for all $x\in(\al_{i_0},x_l]$, $f_a(x)>\al_{i_0}$ for all $x\in[z_a,x_l]$. This, together with (3.4), yields that $f'_a(x)\geq\ld>1$ for all $x\in[z_a,x_l]$. It follows that $f_a(x)>x$ for all $x\in(z_a,x_l]$. So, $f_a(x)<x$ for all $x\in[\al_{i_0},z_a)$, which implies that there exists an integer $k_x\geq1$ such that $f^{k_x}_a(x)<\al_{i_0}$ for any $x\in[\al_{i_0},z_a)$.

For any $a$ satisfying (3.5), one has $f_a(x_r)\geq\al_{i_0+1}$. This, together with $f_a(\al_{i_0+1})=0$ and $f'_a(x)<0$ for all $x\in[x_r,\al_{i_0+1})$, implies that there exists a unique point $u_a$, $u_a\in(x_r,\al_{i_0+1})$, such that $f_a(u_a)=z_a$. On the other hand, there exists a unique point $x_{L,a}$, $x_{L,a}\in(z_a,x_l)$, such that $f_a(x_{L,a})=u_a$, since $f_a(x_l)\geq\al_{i_0+1}$ and $f'_a(x)>0$ for all $x\in(\al_{i_0},x_l)$. From the monotonicity of $f_a(x)$ on $(x_r,\al_{i_0+1})$, it follows that there exists $x_{R,a}$, $x_{R,a}\in(x_r,u_a)$, such that $f_a(x_{R,a})=u_a$.

Thus, for any $a$ satisfying (3.5), there exists a partition of the interval $[\al_{i_0},\al_{i_0+1}]$,
$$[\al_{i_0},\al_{i_0+1}]=[\al_{i_0},z_a)\cup[z_a,x_{L,a}]
\cup(x_{L,a},x_l]\cup(x_l,x_r)\cup[x_r,x_{R,a})\cup[x_{R,a},u_a]
\cup(u_a,\al_{i_0+1}].$$
By the discussions above and the graph of the function $f_a(x)$ on $[\al_{i_0},\al_{i_0+1}]$, one has that for any $x\in[\al_{i_0},\al_{i_0+1}]\backslash([z_a,x_{L,a}]\cup[x_{R,a},u_a])$, there exists an integer $k_x\geq1$ such that
$f^{k_x}_a(x)\not\in[\al_{i_0},\al_{i_0+1}]$.

\indent \bf{Step}\ \rm{2}. The invariant Cantor set is given.

Fix any $a$ satisfying (3.5). Denote
$$I_1:=[z_a,x_{L,a}],\ I_2:=[x_{R,a},u_a],\ I:=I_1\cup I_2.$$
Introduce the notation which is used in the proof of Theorem 4.1 in \cite[Page 30]{Robinson}:
$$I_{i_0,...,i_{n-1}}:=\bigcap^{n-1}_{k=0}f^{-k}_a(I_{i_k})=\{x:\ f^k_a(x)\in I_{i_k}\ \mbox{for}\ 0\leq k\leq n-1\},$$
where $i_k=1$ or $2$, and
$$S_n:=\bigcap^n_{k=0}f^{-k}_a(I)=\bigcap^{n-1}_{k=0}f^{-k}_a(I_1\cup I_2)=\bigcup_{i_0,i_1,...,i_{n-1}=1,2}I_{i_0,i_1,...,i_{n-1}}.$$

It follows from $f_a(x)>\al_{i_0}$ for all $x\in I$ and (3.4) that
$$|f'_a(x)|\geq\ld,\ x\in I.\eqno(3.6)$$

Applying the same method used in the proof of Theorem 4.1 in \cite[Page 30]{Robinson}, one has the following statements, which are similar with those obtained in Lemmas 4.2 and 4.4 in \cite[Page 31]{Robinson}.

\begin{itemize}
\item[{\rm (a)}] For any choice of the labeling with $i_0,...,i_{n-1}\in\{1,2\}$, $I_{i_0,...,i_{n-2}}\cap S_n=I_{i_0,...,i_{n-2},1}\cup I_{i_0,...,i_{n-2},1}$ is the union of two nonempty disjoint closed intervals, which are subsets of $I_{i_0,...,i_{n-2}}$.
\item[{\rm (b)}] For two disjoint choices of the labeling $(i_0,...,i_{n-1})\neq(i'_0,...,i'_{n-1})$, $I_{i_0,...,i_{n-1}}\cap I_{i'_0,...,i'_{n-1}}=\emptyset$, so $S_n$ is the union of $2^n$ disjoint intervals.
\item[{\rm (c)}] The map $f_a$ takes the component $I_{i_0,...,i_{n-1}}$ of $S_n$ homeomorphically onto the component $I_{i_1,...,i_{n-1}}$ of $S_{n-1}$.
\item[{\rm (d)}] The length of any component $I_{i_0,...,i_{n-1}}$ is bounded by $\ld^{-n}(\al_{i_0+1}-\al_{i_0})$.
\end{itemize}

Set
$$\Ld_a:=\bigcap^{\infty}_{n=0}S_n.$$
It is evident that
$$\Ld_a=\{x:\ f^k_a(x)\in[\al_{i_0},\al_{i_0+1}]\ \mbox{for all}\ k\geq0\}.$$
By applying the similar method used in the proof of Theorem 4.1 in \cite[Page 30]{Robinson}, one has that $\Ld_a$ is a Cantor set and a hyperbolic set for $f_a$. From Theorem 5.2 in \cite[Page 38]{Robinson}, it follows that $f_a:\Ld_a\to\Ld_a$ is topologically conjugate to $\si:\sum_2\to\sum_2$.

This completes the whole proof.
\end{proof}

Next, we study the case that $\al_{i_0}=0$ and $f^{(m_{i_0})}_a(0)>0$. For the convenience of our discussions, rewrite $f_a(x)$ as follows:
$$f_a(x)=ax^{d_1}(b-x)^{d_2}h(x),\eqno(3.7)$$
where $d_1=m_{i_0}$, $d_2=m_{i_0+1}$, $b=\al_{i_0+1}$, and $ah(x)>0$ for all $x\in[0,b]$. Fix a constant $0<\ep<1$.

\begin{lemma} \label{lvalue}
Consider the polynomial $f_a(x)$ in the form (3.7).
\begin{itemize}
\item[{\rm (i)}] If $d_1\geq2$, then there exists $\eta>0$ such that $f'_a(x)>0$ and $f''_a(x)>0$ for all $x\in(0,\eta]$;
    $$\frac{f'_a(x)}{f_a(x)}=\frac{g'(x)}{g(x)}\geq\frac{d_1-\ep}{x},\ \ x\in(0,\eta];\ \eqno(3.8)$$
    and there exists a continuous function $u_a:=u(a)$ for $|a|>\eta/|g(\eta)|$ satisfying $0<u_a<\eta$, such that $f_a(u_a)=u_a$, $f_a(x)<x$ for all $x\in(0,u_a)$, and $f_a(x)>x$ for all $x\in(u_a,\eta]$; further, $u_a\to0$ as $|a|\to\infty$.
\item[{\rm (ii)}] If $d_2\geq2$, then there exists $\rho>0$ such that $f'_a(x)<0$ and $f''_a(x)>0$ for all $x\in[\rho,b)$;
    $$-\frac{f'_a(x)}{f_a(x)}=-\frac{g'(x)}{g(x)}\geq\frac{d_2-\ep}{b-x},\ \ x\in[\rho,b);\eqno(3.9)$$
    and there exists a continuous function $v_a:=v(a)$ for $|a|>b/|g(\rho)|$ satisfying $\rho<v_a<b$, such that $f_a(v_a)=b-v_a$, $f_a(x)<b-x$ for all $x\in(v_a,b)$, and $f_a(x)>b-x$ for all $x\in[\rho,v_a)$; further, $v_a\to b$ as $|a|\to\infty$.
\item[{\rm (iii)}] If $d_1>d_2\geq2$, then $u_a>b-v_a$ for sufficiently large $|a|$; if $d_2>d_1\geq2$, then $u_a<b-v_a$ for sufficiently large $|a|$.
\end{itemize}
\end{lemma}

\begin{proof}
Show assertion $\rm(i)$. By Lemma \ref{deriva}, (3.1), (3.2) and (3.3), there exists $\eta>0$ such that $f'_a(x)>0$, $f''_a(x)>0$ and inequality (3.8) holds for all $x\in(0,\eta]$. So, $f_a(x)$ is convex on $(0,\eta]$. Fix the constant $\eta$. Consider the function $H(a,x)=f_a(x)-x$, where $x\in[0,\eta]$, $|a|>\eta/|g(\eta)|$ and $ah(0)>0$. Obviously, $H(a,0)=0$ and $\frac{\partial H}{\partial x}(a,0)=-1<0$. So, for any fixed $|a|>\eta/|g(\eta)|$, there exists $\theta_a>0$ such that $H(a,x)<0$ for all $x\in(0,\theta_a)$. This, together with $H(a,\eta)>0$, implies that there exists $u_0\in(\theta_a,\eta)$ such that $H(a,u_0)=0$ by the intermediate value theorem. Denote $u_a:=\max\{u_0:\  u_0\in(0,\eta)\ \mbox{and}\ H(a,u_0)=0\}$. It is evident that $f_a(u_a)=u_a$ and $f_a(x)>x$ for all $x\in(u_a,\eta]$. Further, This, together with the convexity of $f_a(x)$ on $(0,\eta]$, implies that $f_a(x)<x$ for all $x\in(0,u_a)$. From (3.8), it follows that $\frac{\partial H}{\partial x}(a,u_a)=ag'(u_a)-1\geq(d_1-\ep)-1>0$. Hence, $u_a:=u(a)$ is a continuous function for $|a|>\eta/|g(\eta)|$ by the implicit function theorem. By contradiction, it can be easily concluded that $u_a\to0$ as $|a|\to\infty$. We have now proved that $\rm(i)$ holds.

Show assertion $\rm(ii)$. It follows from Lemma \ref{deriva} (3.1), (3.2) and (3.3) that there exists $\rho>0$ such that $f'_a(x)<0$, $f''_a(x)>0$ and inequality (3.9) holds for all $x\in[\rho,b)$. Thus, $f_a(x)$ is convex on $[\rho,b)$. Fix the constant $\rho$. Consider the function $H(a,x)=f_a(x)+x-b$, where $x\in[\rho,b]$, $|a|>b/|g(\rho)|$, and $ah(0)>0$. It can be easily calculated that $H(a,b)=0$ and $\frac{\partial H}{\partial x}(a,b)=1>0$. So, for any fixed $|a|>b/|g(\rho)|$, there exists $\vartheta_a>0$ such that $H(a,x)<0$ for all $x\in(\vartheta_a,b)$. This, together with $H(a,\rho)>0$, yields that there exists $v_0\in(\rho,b)$ such that $H(a,v_0)=0$ by the intermediate value theorem. Set $v_a:=\min\{v_0:\  v_0\in(0,\rho)\ \mbox{and}\ H(a,v_0)=0\}$. It is evident that $f_a(v_a)=b-v_a$ and $f_a(x)>b-x$ for all $x\in[\rho,v_a)$. This, together with the convexity of $f_a(x)$ on $[\rho,b)$, implies that $f_a(x)<b-x$ for all $x\in(v_a,b)$. By (3.9), $\frac{\partial H}{\partial x}(a,v_a)=ag'(v_a)+1\leq(-d_2+\ep)+1<0$. Hence, it follows from the implicit function theorem that $v_a:=v(a)$ is a continuous function for $|a|>b/|g(\rho)|$. By contradiction, it can be easily shown that $v_a\to b$ as $|a|\to\infty$. Hence, $\rm(ii)$ holds.

Finally, it is to show assertion $\rm(iii)$.

Consider the situation that $d_1>d_2\geq2$. We utilize conclusions and notations in the proof of $\rm(i)$ and $\rm(ii)$ of this lemma. Set
$$N_0:=\max\{\eta/|g(\eta)|,\ b/|g(\rho)|\},\ \ m_1:=\sup_{x,y\in[0,b]}\bigg(\frac{|h(x)|}{|h(y)|}\bigg)^{1/(d_1-d_2)}.$$
Obviously, $0<m_1<+\infty$.

By $\rm(i)$ and $\rm(ii)$ of this lemma, there exists a constant $N_1\geq N_0$ such that if $|a|>N_1$, then
$$\max\{u_a,b-v_a\}<\frac{(d_2-1)b}{d_1+d_2-1},\eqno(3.10)$$
and
$$\frac{b-u_a}{u_a}>m_1,$$
which yields that
$$(b-u_a)^{d_1}u_a^{d_2-1}ah(v_a)>u_a^{d_1-1}(b-u_a)^{d_2}ah(u_a).\eqno(3.11)$$
By $f_a(u_a)=u_a$ and $f_a(v_a)=b-v_a$, one has
$$au_a^{d_1-1}(b-u_a)^{d_2}h(u_a)=av^{d_1}_a(b-v_a)^{d_2-1}h(v_a)=1.\eqno(3.12)$$
Consider the function $G(x)=x^{d_2-1}(b-x)^{d_1}$, $x\in[0,b]$. It is easy to obtain that $G'(x)>0$ for all $x\in(0,(d_2-1)b/(d_1+d_2-1))$.
This, together with (3.10), (3.11) and (3.12), implies that $u_a>b-v_a$. By the same method, we could show that $u_a<b-v_a$ when $d_2>d_1\geq2$ and $|a|$ is sufficiently large. Hence, $\rm(iii)$ holds. The whole proof is complete.
\end{proof}

\begin{theorem} \label{zero}
For the polynomial $f_a(x)=ag(x)$ in the form (3.7), there are four cases to consider: {\rm (i)} $d_1>d_2\geq2$;\ \ {\rm (ii)} $d_1>d_2=1$;\ \ {\rm (iii)} $d_2\geq d_1>1$;\ \ {\rm (iv)} $d_1=d_2=1$. In these cases, if $|a|$ is large enough, then there exists an invariant Cantor set $\Ld_a$ such that $f_a$ is hyperbolic on $\Ld_a$ and $f_a:\Ld_a\to\Ld_a$ is topologically conjugate to $\si:\sum_2\to\sum_2$. In Cases $\rm{(i)}$-$\rm{(iii)}$, any point $x\in(0,b)\backslash\Ld_a$ either escapes from $[0,b]$ or goes to the fixed point $0$ under the iteration of $f_a$; in Case $\rm{(iv)}$, $\Ld_a=\{x:\ f^k_a(x)\in[0,b]\ \mbox{for all}\ k\geq0\}$.
\end{theorem}

\begin{proof}
The whole proof is divided into three parts.

\indent \bf{Part}\ \rm{1}. The properties of the polynomial $f_a(x)$ with $d_1\geq2$ and $d_2\geq2$ are discussed.

For convenience, we still utilize the notations and conclusions in the proof of Lemma \ref{lvalue}. It follows from Lemma \ref{lvalue} that there exist $\eta$ and $\rho$ with $0<\eta<\rho<b$, such that $\rm(i)$ and $\rm(ii)$ of Lemma \ref{lvalue} hold. Hence, for any $|a|>N_0=\max\{\eta/|g(\eta)|,\ b/|g(\rho)|\}$, one has
$$|f'_a(x)|\geq d_1-\ep,\ x\in[u_a,\eta];\eqno(3.15)$$
$$|f'_a(x)|\geq d_2-\ep,\ x\in[\rho,v_a].\eqno(3.16)$$
Denote
$$m_0:=\inf_{x\in[\eta,\rho]}|g(x)|.$$
Obviously, $0<m_0\leq\min\{|g(\eta)|,\ |g(\rho)|\}.$

In the following discussions, fix any
$$|a|\geq \frac{b}{m_0}\ \mbox{and}\ ag^{(d_1)}(0)>0.\eqno(3.17)$$
Thus, for any fixed $a$ satisfying (3.17), there exist two points $x_{l,a}$ and $x_{r,a}$ with $0<x_{l,a}\leq\eta<\rho\leq x_{r,a}<b$, such that $f_a(x_{l,a})=f_a(x_{r,a})=b$. It is evident that $x_{l,a}\to 0$ and $x_{r,a}\to b$ as $|a|\to+\infty$.

It follows from $f'_a(x)<0$ for all $x\in[\rho,b)$ that there exists a unique point $w_a$, $w_a\in(x_{r,a},b)$, such that $f_a(w_a)=u_a$. On the other hand, there exists a unique point $x_{L,a}$, $x_{L,a}\in(u_a,x_{l,a})$, such that $f_a(x_{L,a})=w_a$, since $f_a(x_{l,a})=b$ and $f'_a(x)>0$ for all $x\in(0,\eta]$. From the monotonicity of $f_a(x)$ on $(\rho,b)$ and $f_a(x_{r,a})=b$, it follows that there exists a unique point $x_{R,a}$, $x_{R,a}\in(x_{r,a},w_a)$, such that $f_a(x_{R,a})=w_a$.

Hence, for any $a$ satisfying (3.17), there exists a partition of the interval $(0,b)$,
$$(0,b)=(0,u_a)\cup[u_a,x_{L,a}]
\cup(x_{L,a},\eta]\cup(\eta,\rho)\cup[\rho,x_{R,a})\cup[x_{R,a},w_a]
\cup(w_a,b).$$
Since $f_a(x)<x$ for all $x\in(0,u_a)$, $f_a^k(x)\to0$ as $k\to+\infty$ for any $x\in(0,u_a)$. So, from the graph of the function $f_a(x)$ on $[0,b]$, it follows that for any $x\in(0,b)\backslash([u_a,x_{L,a}]\cup[x_{R,a},w_a])$, either there exists an integer $k_x\geq1$ such that
$f^{k_x}_a(x)\not\in(0,b)$ or $f_a^k(x)\to0$ as $k\to+\infty$. Denote
$$I_1:=[u_a,x_{L,a}],\ I_2:=[x_{R,a},w_a],\ I:=I_1\cup I_2,\ S_n:=\bigcap^n_{k=0}f^{-k}_a(I),\ \Ld_a:=\bigcap^{\infty}_{n=0} S_n.$$
Obviously,
$$f_a(I_1)=f_a(I_2)\supset I_1\cup I_2.$$
It is easy to conclude that $f_a:\Ld_a\to\Ld_a$ is topologically semiconjugate to $\si:\sum_2\to\sum_2$.

Consider Case $\rm{(i)}$. By $\rm(iii)$ of Lemma \ref{lvalue}, there exists $N_1$ such that if $|a|>N_1$, then $u_a>b-v_a$, where $N_1$ is given in the proof of $\rm(iii)$ of Lemma \ref{lvalue}. Hence, $w_a<v_a$ if $|a|>N_1$. This, together with (3.15) and (3.16), implies that when $|a|>N_1$, one has
$$|f'_a(x)|\geq d_1-\ep>1,\ x\in I_1;\eqno(3.18)$$
$$|f'_a(x)|\geq d_2-\ep>1,\ x\in I_2.\eqno(3.19)$$
Applying the same method used in Step $\rm {2}$ of the proof of Theorem \ref{2positive}, one could show the assertion in Case $\rm{(i)}$. By using the similar method in Case $\rm(i)$ and Lemma \ref{singleroot}, one could prove the assertions in Cases $\rm(ii)$ and $\rm(iv)$.

However, when $d_2\geq d_1$, the value of $|f'_a|$ on $I_2$ does not become
very large as $|a|\to+\infty$, which will be explained below.

\indent \bf{Part}\ \rm{2}. Estimate the value of $f'_a(w_a)$ as $|a|\to+\infty$ when $d_1\geq2$ and $d_2\geq2$.

Fix any $a$ satisfying (3.17) in the following discussions. Since
$$f_a(u_a)=u^{d_1}_a(b-u_a)^{d_2}ah(u_a)=u_a,\ \ f_a(w_a)=w^{d_1}_a(b-w_a)^{d_2}ah(w_a)=u_a,$$
$$f'_a(x)=d_1x^{d_1-1}(b-x)^{d_2}ah(x)-d_2x^{d_1}(b-x)^{d_2-1}ah(x)+x^{d_1}(b-x)^{d_2}ah'(x),$$
one has
$$f'_a(w_a)=\frac{d_1u_a}{w_a}-\frac{d_2u_a}{b-w_a}+u_a\frac{h'(w_a)}{h(w_a)}.\eqno(3.20)$$
Rewrite $f_a(x)$ as follows:
$$f_a(x)=ab^{d_2}h(0)x^{d_1}+a\frac{g^{(d_1+1)}(\tau x)}{(d_1+1)!}x^{d_1+1},\ \ x\in[0,\eta], \ \ 0<\tau<1;\eqno(3.21)$$
$$f_a(x)=ab^{d_1}h(b)(b-x)^{d_2}+a\frac{g^{(d_2+1)}(b+\nu(x-b))}{(d_2+1)!}(x-b)^{d_2+1},\ \ x\in[\rho,b],\ \ 0<\nu<1.\eqno(3.22)$$
Set
$$\alpha_1:=b^{d_2}h(0),\  \beta_1(x):=\frac{g^{(d_1+1)}(\tau x)}{(d_1+1)!},\ \alpha_2:=b^{d_1}h(b),\  \beta_2(x):=\frac{g^{(d_2+1)}(b+\nu(x-b))}{(d_2+1)!}.$$
By (3.21) and $f_a(u_a)=u_a$, one has
$$u_a=\bigg(\frac{1}{a\alpha_1+a\beta_1(u_a)u_a}\bigg)^{1/(d_1-1)}.\eqno(3.23)$$
From (3.22) and $f_a(w_a)=u_a$, it follows that
$$b-w_a=\bigg(\frac{u_a}{a\alpha_2+(-1)^{d_2+1}a\beta_2(w_a)(b-w_a)}\bigg)^{1/d_2}.\eqno(3.24)$$
By (3.23) and (3.24), one has
$$\frac{u_a}{b-w_a}=|a|^{\frac{d_1-d_2}{(d_1-1)d_2}}
(|\alpha_1|+|\beta_1(u_a)|u_a)^{\frac{1-d_2}{(d_1-1)d_2}}(|\alpha_2|+(-1)^{d_2+1}|\beta_2(w_a)|(b-w_a))^{\frac{1}{d_2}}.\eqno(3.25)$$

It is evident that $u_a\to0$ and $w_a\to b$ as $|a|\to+\infty$, and $\sup_{x\in[0,\eta]}|\beta_1(x)|$, $\sup_{x\in[\rho,b]}|\beta_2(x)|$, and $\sup_{x\in[0,b]}\frac{|h'(x)|}{|h(x)|}$ are finite. By (3.20), when $|a|$ is large enough, one has
$$f'_a(w_a)=-d_2|a|^{\frac{d_1-d_2}{(d_1-1)d_2}}
(|\alpha_1|^{\frac{1-d_2}{(d_1-1)d_2}}|\alpha_2|^{\frac{1}{d_2}}+o(1))+o(1),\eqno(3.26)$$
where $o(1)\to0$ as $|a|\to+\infty$.

By $f''_a(x)>0$ and $f'_a(x)<0$, $x\in[\rho,b)$, one has that $|f'(w_a)|=\inf_{x\in[x_{R,a},w_a]}|f'(x)|$. This, together with (3.26), yields that $f'_a(w_a)\to0$ as $|a|\to\infty$, when $d_2>d_1\geq2$. This means that the method used in the discussions of Case $\rm(i)$ of this theorem could not be applied to study this case. Hence, we utilize the tools in complex dynamics.

\indent \bf{Part}\ \rm{3}. The dynamical behaviors of the complex polynomial map $f_a(z)=ag(z)=az^{d_1}(b-z)^{d_2}h(z)$ on $\overline{\mathbf{C}}$ are studied, where $a$ and $z$ are complex numbers and $d_1\geq2$.

Now, it is to show that $f_a$ is hyperbolic on its Julia set when $|a|$ is large enough. For a constant $N>0$, denote $D(\infty,N):=\{z\in\mathbf{C},\ |z|>N\}$. It is evident that for the complex polynomial $g(z)$ on $\overline{\mathbf{C}}$, there exists a constant $K>0$ such that $g(D(\infty,K))\subset D(\infty,K)$, and $D(\infty,K)$ is contained in the attracting neighborhood of $\infty$ for $g$. Fix this constant $K$. So, for any $|a|\geq1$, $f_a(D(\infty,K))\subset D(\infty,K)$, and $D(\infty,K)$ is contained in the attracting neighborhood of $\infty$ for $f_a$. The critical set $E:=\{z\in\mathbf{C},\ g'(z)=0\}$ of $g$ is finite. And, $E=E_1\cup E_2$, where $E_1=\{z\in E,\ g(z)\neq0\}$, and $E_2=\{z\in E,\ g(z)=0\}$. Since $d_1\geq2$, $f_a(0)=0$, and $f_a(b)=0$, one has that $0$ is a superattracting fixed point of $g$, $0\in E_2$, and $E_1\neq\emptyset$. Obviously, $E_2\subset A(0,f_a)$, and, if
$$|a|>\max\bigg\{1,\ \max_{z\in E_1}\frac{K}{|g(z)|}\bigg\},\eqno(3.27)$$
then $f_a(E_1)\subset D(\infty,K)\subset A(\infty,f_a)$. Hence, one has that $f_a$ is hyperbolic on the Julia set $\emph{J}$ for any $a$ satisfying (3.27).

It is to show that $\Ld_a$ is contained in the Julia set for large $|a|$. It can be easily obtained that $A(0,f_a)$ and $A(\infty,f_a)$ are two disjoint completely invariant components of the Fatou set. So, $\partial A(0,f_a)=\partial A(\infty,f_a)=\emph{J}$. And, there are at most two completely invariant components of the Fatou set. So, the Fatou set coincides with $A(0,f_a)\cup A(\infty,f_a)$. When $f_a$ is restricted on the real line, one has that if $a$ satisfies both (3.17) and (3.27), then $\Ld_a\not\subset A(0,f_a)\cup A(\infty,f_a)$, which implies that $\Ld_a\subset \emph{J}$.

Hence, $\Ld_a$ is a hyperbolic invariant set for $f_a$ for $d_1\geq2$ and sufficiently large $|a|$. Thus, for any $|a|$ satisfying (3.17) and (3.27), there exist two constants $C_a>0$ and $\ld_a>1$ such that $|(f^k_a)'(x)|\geq C_a\ld^k_a$ for all $k\geq0$ and $x\in\Ld_a$. This, together with $f_a:\Ld_a\to\Ld_a$ is topologically semiconjugate to $\si:\sum_2\to\sum_2$, yields that $f_a:\Ld_a\to\Ld_a$ is topologically conjugate to $\si:\sum_2\to\sum_2$ for any $|a|$ satisfying (3.17) and (3.27), and $\Ld_a$ is a Cantor set. Therefore, assertion in Case $\rm(iii)$ holds.

The whole proof is complete.
\end{proof}

\begin{remark}
It is a meaningful question to find an elementary method to prove the assertion in Case $\rm(iii)$ of Theorem \ref{zero} without the help of the tools in complex dynamics, since the method in complex dynamics requires lots of preparation.
\end{remark}

Finally, we study the situation that $d_2>d_1=1$.

\begin{theorem} \label{nonhyp}
Consider the polynomial $f_a(x)$ in the form (3.7). Suppose that $d_2>d_1=1$. Then there exists a constant $N>0$ such that for any $|a|>N$, there exists an invariant set $\Ld_a=\{x:\ f^k_a(x)\in[0,b]\ \mbox{for all}\ k\geq0\}$ such that $f_a:\Ld_a\to\Ld_a$ is topologically semiconjugate to $\si:\sum_2\to\sum_2$, but $f_a$ is not hyperbolic on $\Ld_a$.
\end{theorem}

\begin{proof}
By Lemma \ref{singleroot} and the properties of the function $f_a$ on $[0,b]$, it can be easily concluded that there exists a constant $N>0$, such that for any $|a|>N$, there exist two intervals $I_1=[0,\eta]$ and $I_2=[\rho,b]$, such that $f_a(I_1)=f_a(I_2)=[0,b]$, $f'_a(x)>1$ for all $x\in I_1$, and $f_a(x)>b$ for any $x\in(\eta,\rho)$, where $0<\eta<\rho<b$. Set
$$I:=I_1\cup I_2,\ S_n:=\bigcap^n_{k=0}f^{-k}_a(I),\ \Ld_a:=\bigcap^{\infty}_{n=0}S_n.$$
By simple discussions, one has $f_a:\Ld_a\to\Ld_a$ is topologically semiconjugate to $\si:\sum_2\to\sum_2$.

Now, it is to show that $\Ld_a$ is not a hyperbolic invariant set for $f_a$.
Since $b\in\Ld_a$ and $f'_a(b)=0$, one has that $|(f^k_a)'(b)|=0$ for any positive integer $k$. Hence, $\Ld_a$ is not a hyperbolic invariant set for $f_a$ by Lemma \ref{hyp}. The proof is complete.
\end{proof}

\begin{remark}
For the polynomial $f_a(x)$ in the form (3.7) with $d_2>d_1=1$, we conjecture that there does not exist a hyperbolic invariant set $\Ld'_a\varsubsetneqq\Ld_a$ for $f_a$, such that $f_a:\Ld'_a\to\Ld'_a$ is topologically conjugate to $\si:\sum_2\to\sum_2$ for $|a|>N$, where $N$ is specified in Theorem \ref{nonhyp}.
\end{remark}

Finally, we consider the situation that $\al_{i_0}\geq0$ and $f^{(m_{i_0})}_a(\al_{i_0})<0$. There are four different cases to consider:
$\rm(a)$ $m_{i_0}$ is odd and $f^{(m_{i_0+1})}_a(\al_{i_0+1})>0$; $\rm(b)$ $m_{i_0}$ is odd and $f^{(m_{i_0+1})}_a(\al_{i_0+1})<0$;  $\rm(c)$ $m_{i_0}$ is even and $f^{(m_{i_0+1})}_a(\al_{i_0+1})>0$; $\rm(d)$ $m_{i_0}$ is even and $f^{(m_{i_0+1})}_a(\al_{i_0+1})<0$.

\begin{theorem} \label{2positiveneg}
Suppose that $\al_{i_0}>0$ and $f^{(m_{i_0})}_a(\al_{i_0})<0$. In Case $\rm(a)$, there exists an invariant set $\Ld_a$ contained in $(\al_{i_0-1},\al_{i_0})\cup(\al_{i_0+1},\al_{i_0+2})$ such that $f_a$ is hyperbolic on $\Ld_a$ and $f_a:\Ld_a\to\Ld_a$ is topologically conjugate to $\si:\sum_2\to\sum_2$ for sufficiently large $|a|$, where $\al_{i_0-1}=-\infty$ if $i_0=1$, $\al_{i_0+2}=+\infty$ if $i_0=r-1$.
\end{theorem}

\begin{proof}
By a modification of the method used in the proof of Case $\rm(1)$ in Theorem 4.1, one can verify the conclusions of this theorem. So, we only give the idea.  Fix a constant $\ld>1$. When $|a|$ is large enough, we could find two disjoint compact intervals $I_1$ and $I_2$, $I_1\subset (\al_{i_0-1},\al_{i_0})$, $I_2\subset (\al_{i_0+1},\al_{i_0+2})$, such that $f_a(I_1)=f_a(I_2)\supset I_1\cup I_2$, and $|f'_a(x)|\geq\ld$ for all $x\in I_1\cup I_2$. Repeat the discussions in Step 2 of Theorem \ref{2positive}, one could show the conclusions of this theorem.
\end{proof}

However, in the situation that $\al_{i_0}=0$ and $f^{(m_{i_0})}_a(\al_{i_0})<0$, it is easy to obtain the following result.

\begin{theorem} \label{zeropositive}
In Case $\rm(a)$, for sufficiently large $|a|$, there exists an invariant set $\Ld_a\subset (\al_{i_0-1},\al_{i_0+2})$ for $f_a$, on which $f_a$ is topologically semiconjugate to $\si:\sum_2\to\sum_2$ ($\Ld_a$ might not be hyperbolic); further, if $m_{i_0}=m_{i_0+1}=1$, $\Ld_a$ is a hyperbolic invariant set for $f_a$. In Case $\rm(b)$, for sufficiently large $|a|$, there exists an invariant set $\Ld_a\subset(\al_{i_0-1},\al_{i_0+1})$ for $f_a$ such that $f_a:\Ld_a\to\Ld_a$ is topologically semiconjugate to $\si_A:\sum_A\to\sum_A$, where
\[A=\left(
  \begin{array}{cc}
    1 & 1 \\
    1 & 0 \\
  \end{array}
\right),
\]
and it is an eventually positive transition matrix.
\end{theorem}
\bigskip

\section{Polynomial with one positive and one negative zeros}

In this section, two classes of real-coefficient polynomial maps are investigated, where one class has one positive and one negative zeros, and the other has only one real zero.

Consider the real-coefficient polynomial $g(x)$ introduced in Section 3. In this section, suppose that there exists $i_0$, $1\leq i_0<r$, such that $\al_{i_0}<0<\al_{i_0+1}$, we study the polynomial map $f_a(x)=ag(x)$,
where $a\in \mathbf{R}$ is a parameter.

There are eight different cases to consider:
\begin{itemize}
\item[{\rm (1)}] $f^{(m_{i_0})}_a(\al_{i_0})>0$, $f^{(m_{i_0+1})}_a(\al_{i_0+1})>0$, $m_{i_0}$ is even;
\item[{\rm (2)}] $f^{(m_{i_0})}_a(\al_{i_0})>0$, $f^{(m_{i_0+1})}_a(\al_{i_0+1})<0$, $m_{i_0}$ is even;
\item[{\rm (3)}] $f^{(m_{i_0})}_a(\al_{i_0})>0$, $f^{(m_{i_0+1})}_a(\al_{i_0+1})>0$, $m_{i_0}$ is odd;
\item[{\rm (4)}] $f^{(m_{i_0})}_a(\al_{i_0})>0$, $f^{(m_{i_0+1})}_a(\al_{i_0+1})<0$, $m_{i_0}$ is odd;
\item[{\rm (5)}] $f^{(m_{i_0})}_a(\al_{i_0})<0$, $f^{(m_{i_0+1})}_a(\al_{i_0+1})>0$, $m_{i_0}$ is odd;
\item[{\rm (6)}] $f^{(m_{i_0})}_a(\al_{i_0})<0$, $f^{(m_{i_0+1})}_a(\al_{i_0+1})<0$, $m_{i_0}$ is odd;
\item[{\rm (7)}] $f^{(m_{i_0})}_a(\al_{i_0})<0$, $f^{(m_{i_0+1})}_a(\al_{i_0+1})>0$, $m_{i_0}$ is even;
\item[{\rm (8)}] $f^{(m_{i_0})}_a(\al_{i_0})<0$, $f^{(m_{i_0+1})}_a(\al_{i_0+1})<0$, $m_{i_0}$ is even.
\end{itemize}

\begin{theorem} \label{diffroot}
In Cases $\rm(2)$ and $\rm(6)$, for any sufficiently large $|a|$, there exists a hyperbolic invariant set $\Ld_a$ for $f_a$ such that $f_a:\Ld_a\to\Ld_a$ is topologically conjugate to $\si_A:\sum_A\to\sum_A$, where
\[A=\left(
  \begin{array}{ccc}
    0 & 0 & 1 \\
    0 & 0 & 1 \\
    1 & 1 & 1 \\
  \end{array}
\right).
\]
In other cases, for any sufficiently large $|a|$, there exists a hyperbolic invariant set $\Ld_a$ for $f_a$ such that $f_a:\Ld_a\to\Ld_a$ is topologically conjugate to $\si:\sum_2\to\sum_2$.
\end{theorem}

\begin{remark}
It can be directly calculated that $A^2$ is positive. So, $A$ is an eventually positive transition matrix, and $\si_A$ is topologically mixing on $\sum_A$ \cite[Page 77, Proposition 2.9]{Robinson}.
\end{remark}

\begin{proof}
We only give the proofs of Cases $\rm(1)$ and $\rm(2)$, other cases can be studied by using the similar method.

Consider Case $\rm(1)$. By assumptions, one has that $f_a(x)>0$ for all $x\in (\al_{i_0-1},\al_{i_0})\cup(\al_{i_0},\al_{i_0+1})\cup(\al_{i_0+1},\al_{i_0+2})$, where $\al_{i_0-1}=-\infty$ if $i_0=1$, $\al_{i_0+2}=+\infty$ if $i_0=r-1$. Fix a constant $\ld>1$. By (3.2), there exist $x_1$ and $x'_1$, $\al_{i_0+1}/2\leq x_1<\al_{i_0+1}<x'_1$, such that $f'_a(x)<0$ for all $x\in[x_1,\al_{i_0+1})$, $f'_a(x)>0$ for all $x\in(\al_{i_0+1},x'_1]$, and
$$\bigg|\frac{ f_a'(x)}{
f_a(x)}\bigg|=\bigg|\frac{g'(x)}{g(x)}\bigg|\geq\frac{2\ld}{\al_{i_0+1}},\ \ x\in[x_1,\al_{i_0+1})\cup(\al_{i_0+1},x'_1].         \eqno{(4.1)}$$

In the following discussions, fix any
$$|a|\geq\frac{x'_1}{\min\{|g(x_1)|,|g(x'_1)|\}}.$$
Hence, $f_a(x'_1)\geq x'_1$, which together with $f_a(\al_{i_0+1})=0$, implies that there exists $y\in(\al_{i_0+1},x'_1]$ such that $f_a(y)=y$. Denote $x_{R,a}:=\max\{y:\ f_a(y)=y,\ y\in(\al_{i_0+1},x'_1]\}$. So, $f_a(x_{R,a})=x_{R,a}$. Since $f_a(x_1)\geq x'_1$, $f_a(\al_{i_0+1})=0$, and $f'_a(x)<0$ for all $x\in[x_1,\al_{i_0+1})$, there exist unique two points $x_{L,a}$ and $x_{l,a}$, $x_1\leq x_{L,a}<x_{l,a}<\al_{i_0+1}$, such that $f_a(x_{L,a})=x_{R,a}$ and $f_a(x_{l,a})=x_{L,a}$. And, there exists $x_{r,a}$, $x_{r,a}\in(\al_{i_0+1},x_{R,a})$, such that $f_a(x_{r,a})=x_{L,a}$. Set
$$I_1:=[x_{L,a},x_{l,a}],\ I_2:=[x_{r,a},x_{R,a}].$$
Thus,
$$f_a(I_1)=f_a(I_2)\supset I_1\cup I_2.$$
By (4,1), one has
$$|f'_a(x)|\geq\ld,\ x\in I_1\cup I_2.$$
By applying the method used in Step 2 of Theorem \ref{2positive}, one could obtain the conclusions in this case.

Now, it is to consider Case $\rm(2)$. It follows from the hypothesis that $f_a(x)>0$ for all $x\in (\al_{i_0-1},\al_{i_0})\cup(\al_{i_0},\al_{i_0+1})$ and $f_a(x)<0$ for all $x\in(\al_{i_0+1},\al_{i_0+2})$, where $\al_{i_0-1}=-\infty$ if $i_0=1$, $\al_{i_0+2}=+\infty$ if $i_0=r-1$.
Choose a constant $\ld>1$. It follows from (3.2) that there exist $y_1$, $y_2$, $y_3$ and $y_4$, such that $f'_a(x)<0$ for all $x\in[y_1,\al_{i_0})\cup[y_3,\al_{i_0+1})\cup(\al_{i_0+1},y_4]$, $f'_a(x)>0$ for all $x\in(\al_{i_0},y_2]$, and
$$\bigg|\frac{ f_a'(x)}{
f_a(x)}\bigg|=\bigg|\frac{g'(x)}{g(x)}\bigg|\geq\frac{2\ld}{\min\{\al_{i_0+1},|\al_{i_0}|\}},\ \ x\in[y_1,\al_{i_0})\cup(\al_{i_0},y_2]\cup[y_3,\al_{i_0+1})\cup(\al_{i_0+1},y_4],         \eqno{(4.2)}$$
where $y_1<\al_{i_0}<y_2\leq \al_{i_0}/2<\al_{i_0+1}/2\leq y_3<\al_{i_0+1}<y_4$.

In the following discussions , fix any
$$|a|\geq\frac{\max\{|y_1|,y_4\}}{\min_{1\leq i\leq 4}\{|g(y_i)|\}}.$$
By intermediate value theorem, there exist $y_{1,a}$ and $y_{2,a}$, $y_1\leq y_{1,a}<\al_{i_0}<y_{2,a}<y_2$, such that $f_a(y_{1,a})=f_a(y_{2,a})=y_3$. Denote
$$I_1:=[y_1,y_{1,a}],\ I_2:=[y_{2,a},y_2],\ I_3:=[y_3,y_4].$$
Hence,
$$f_a(I_1)\supset I_3,\ f_a(I_2)\supset I_3,\ f_a(I_3)\supset I_1\cup I_2\cup I_3.$$
For any $\be=(b_0,b_1,b_2,...)\in\sum_A$, set
$$I_{b_0,...,b_{n-1}}:=\bigcap^{n-1}_{k=0}f^{-k}_a(I_{b_k})=\{x:\ f^k_a(x)\in I_{b_k}, \ 0\leq k\leq n-1\},$$
$$I_{\be}:=\bigcap^{+\infty}_{n=0}I_{b_0,...,b_n}.$$
It can be concluded that $I_{\be}$ is either a nondegenerate compact interval or a singlton set, since $f_a$ is monotone on $I_i$, $1\leq i\leq 3$. For any $\be=(b_0,b_1,b_2,...)\in\sum_A$, one has that
$$f(I_{\be})=\bigcap^{\infty}_{n=0}f(I_{b_0\cdots b_n})=\bigcap^{\infty}_{n=1}I_{b_1\cdots b_n}=I_{\si_A(\be)}.\eqno(4,3)$$
By induction, one has that $I_{c_1\cdots c_k}\bigcap I_{d_1\cdots d_k}=\emptyset$ for any two different allowable words $w_1=(c_1,...,c_k)$ and $w_2=(d_1,...,d_k)$
for $A$. Consequently, for any
$\be,\ga\in\sum_A$ with $\be\neq\ga$, one has
$$I_{\be}\cap I_{\ga}=\emptyset.$$
Set
$$\Ld_a:=\bigcup_{\be\in\sum_A} {I_{\be}}.$$
It follows from (4.3) that $f_a(\Ld_a)=\Ld_a$ since $A$ is a transition matrix. Hence, $f_a(x)\in\Ld_a$ for any $x\in\Ld_a$, which yields that $\al_{i_0+1}\not\in\Ld_a$, and $|f_a(x)|\geq\min\{|\al_{i_0}|/2,\al_{i_0+1}/2\}$ for any $x\in\Ld_a$. This, together with (4.2), implies that
$$|f'_a(x)|\geq\ld,\ x\in\Ld_a.\eqno(4.4)$$

Now, it is to show that $I_{\be}$ is a singlton set for any $\be\in\sum_A$.
By contradiction, assume that $I_{\be}$ a nondegenerate compact interval. By (4.3), (4.4) and $f_a$ is monotone on $I_i$, $1\leq i\leq3$, one has
$$|I_{\si^k_A(\be)}|=|f^k_a(I_{\be})|\geq\ld^k|I_{\be}|,\ k\geq1,$$
which yields that $|I_{\si^k_A(\be)}|\to+\infty$ as $k\to+\infty$. It is a contradiction. Thus, $I_{\be}$ is a singlton set.

Therefore, $\Ld_a$ is a hyperbolic invariant set for $f_a$ and $f_a:\Ld_a\to\Ld_a$ is topologically conjugate to $\si_A:\sum_A\to\sum_A$ in Case $\rm(2)$.

One could apply the method used in the discussions of Case $\rm(1)$ to study Case $\rm(3)$. However, in the case of $\rm(4)$, for sufficiently large $|a|$, one could find to two disjoint compact intervals $I_1$ and $I_2$ such that $I_1\cup I_2\subset f_a(I_1), f_a(I_2)$, where $\al_{i_0}\in I_1\subset (\al_{i_0-1},0)$ and $\al_{i_0+1}\in I_2\subset (0,\al_{i_0+2})$. By using the method used in the study of Case $\rm(2)$, one could prove that the invariant set contained in $I_1\cup I_2$ is a hyperbolic invariant set for $f_a$ on which $f_a$ is topologically conjugate to $\si:\sum_2\to\sum_2$. Cases $\rm(5)$-$\rm(8)$ can be studied with similar discussions.

This completes the whole proof.
\end{proof}

By applying the method used in the discussions of Case $\rm(1)$ in Theorem \ref{diffroot}, one could easily obtain the following result.

\begin{theorem} \label{single}
For the polynomial map $f_a(x)=ag(x)$, suppose that $g$ has only one real zero $\al$. If $\al>0$ and $f_a(x)\geq0$ for all $x\in\mathbf{R}$, then for sufficiently large $|a|$, there exists a hyperbolic invariant set $\Ld_a\subset(0,\al)\cup(\al,+\infty)$ for $f_a$ and $f_a:\Ld_a\to\Ld_a$ is topologically conjugate to $\si:\sum_2\to\sum_2$. If $\al<0$ and $f_a(x)\leq0$ for all $x\in\mathbf{R}$, then for sufficiently large $|a|$, there exists a hyperbolic invariant set $\Ld_a\subset(-\infty,\al)\cup(\al,0)$ for $f_a$ and $f_a:\Ld_a\to\Ld_a$ is topologically conjugate to $\si:\sum_2\to\sum_2$.
\end{theorem}


\begin{thebibliography}{30}
%\baselineskip=15pt
\bibitem{Aulbach} B. Aulbach, B. Kieninger, An elementary proof for hyperbolicity and chaos of the logistic maps, J. Diff. Equ. Appl. 10 (2004) 1243--1250.

\bibitem{CarlesonGamelin} L. Carleson, T. W. Gamelin, Complex Dynamics, Springer-Verlag, New York, 1993.

\bibitem{Glendinning} P. Glendinning, Hyperbolicity of the invariant set for the logistic map
with $\mu>4$, Nonlinear Analysis 47 (2001) 3323--3332.

\bibitem{Guckenheimer} J. Guckenheimer, Sensitive dependence to initial conditions for one dimensional maps, Commun. Math. Phys. 70 (1979) 133--160.

\bibitem{Henry} B. R. Henry, Escape from the unit interval under the transformation $x\to\ld x(1-x)$, Proc. Am. Math. Soc. 41 (1973) 146--150.

\bibitem{Kraft} R. L. Kraft, Chaos, Cantor sets, and hyperbolicity for the logistic maps, Am. Math. Mon. 106 (1999) 400--408.

\bibitem{Misiurewicz} M. Misiurewicz, Absolutely continuous measures for certain maps of an interval, Publ. Math. I.H.E.S. 53 (1981) 17--51.

\bibitem{Robinson} C. Robinson, Dynamical Systems:
Stability, Symbolic Dynamics and Chaos, CRC Press, Boca Raton, 1995.

\bibitem{Strien} S. van Strien, On the bifurcation creating horseshoes, in: Rand and Young (Eds.), Dynamical Systems and Turbulence, Lecture Notes in Math., vol. 898, Springer-Verlag, New York, Heidelberg, Berlin, 1981, pp. 316--351.

\end{thebibliography}
\end{document}